\documentclass[11pt]{article}
\usepackage{amssymb}
\usepackage{amsfonts}
\usepackage{amscd}

\usepackage[english]{babel}

\usepackage{amsmath,latexsym,enumerate,amsthm,amssymb}

\usepackage{calc}

\newcommand{\be}{\begin{enumerate}}
\newcommand{\ee}{\end{enumerate}}
\let\ds=\displaystyle

\def\N{{\mathbb N}} \def\Z{{\mathbb Z}}

 \def\Q{{\mathbb Q}}

\def\R{{\mathbb R}} \def\C{{\mathbb C}}

\def\P{{\mathbb P}}
\def\Sp{{\mathbb S}}

\def\T{{\cal T}}

\def\s{{\bf s}}

\def\m{{\bf m}}
\def\i{{\bf i}}

\def\z{{\bf z}}
\def\e{{\bf e}}
\def\x{{\bf x}}
\def\y{{\bf y}}
\def\a{{\bf a}}
\def\cb{{\bf c}}
\def\r{{\bf r}}
\def\u{{\bf u}}
\def\v{{\bf v}}

\def\b{{\bf b}}
\def\m{{\bf m}}
\def\X{{\bf X}}

\newcommand{\zerob}{\boldsymbol{0}}
\newcommand{\unb}{\boldsymbol{1}}

\newcommand{\Ab}{\mathbf{A}}

\newcommand{\alphab}{\boldsymbol{\alpha}}
\newcommand{\betab}{\boldsymbol{\beta}}
\newcommand{\gammab}{\boldsymbol{\gamma}}

\newcommand{\nub}{\boldsymbol{\nu}}
\newcommand{\sigmab}{\boldsymbol{\sigma}}
\newcommand{\taub}{\boldsymbol{\tau}}
\newcommand{\lambdab}{\boldsymbol{\lambda}}

\newcommand{\etab}{\boldsymbol{\eta}}

\def\cqfd{ $\diamondsuit $ }

\def\eps{{ \varepsilon }}

\catcode`\ç=13
\defç{\c{c}}
\catcode`\é=13
\defé{\'e}
\catcode`\à=13
\defà{\`a}
\catcode`\è=13
\defè{\`e}
\catcode`\â=13
\defâ{\^a}
\catcode`\ù=13
\defù{\`u}
\catcode`\ê=13
\defê{\^e}
\catcode`\î=13
\defî{\^\i}
\catcode`\ô=13
\defô{\^o}

\newcommand{\la}{\langle}

\newcommand{\ra}{\rangle}

\newtheorem{Lemme}{Lemma}

\newtheorem{Corollaire}{Corollary}

\newtheorem{Proposition}{Proposition}

\newtheorem{Theoreme}{Theorem}

\newtheorem{Definition}{Definition}

\begin{document}
\title{ Manin's conjecture on toric varieties with different heights.} 
\author{Driss Essouabri\footnote{Universit\'e de Caen,
 UFR des Sciences, Campus 2,
 Laboratoire de Math. Nicolas Oresme, CNRS UMR 6139,
 Bd. Mal Juin, B.P. 5186, 14032 Caen, France.
Email : essoua@math.unicaen.fr}\hskip 0.4cm
}
\date{}
\maketitle 

\noindent
{\small {{\bf {Abstract.}}
In this paper I verify Manin's conjecture for a class of rational projective 
toric varieties with a large class of heights other than  the usual 
one that comes from the standard metric on projective space.}\par 
\vskip .2 in
\setcounter{tocdepth}{2}
 \section{Introduction and detailed description of the problem} 

Let $X$ be a projective algebraic variety over $\Q$ and $U$ a Zariski open subset  of $X.$  We can study the
density of rational points of a projective embedding of  $U$  by first choosing a suitable line 
bundle ${\mathcal {L}}$\footnote {technically, one should have that 
the  class of  $ \mathcal{L} $ in the Picard group of $X$ is contained in
the interior of the cone of effective divisors} that determines a 
projective embedding $\phi: U(\Q) \hookrightarrow \P^{n}(\Q)$ for some
$n$. There are then several ways to measure the density of $\phi
\,U(\Q)$. First, we choose a family of metrics  
  ${\it v}=\left(v_p\right)_{p  \le \infty}$ whose Euler product     
 $H_{{\it v}}=\prod_{p \le \infty} v_p$
specifies  a height $H_{v}$ on $\phi \, U(\Q)$ in the evident way. It suffices to set
$H_{{\mathcal {L}}, {\it v}}(M):=H_{{\it v}}\left(\phi(M)\right)$ for any point $M \in U(\Q)$.\\ 
By definition, the density of the rational points  of $U(\Q)$ with respect to
$H_{{\mathcal {L}}, {\it v}}$ is the function  
$$B \to N_{{\mathcal {L}}, {\it v}}(U, B):=\{M \in U(\Q) ~|~H_{{\mathcal{L}}, {\it v}}(M)\leq B\}.$$
  Manin's conjecture concerns  the asymptotic behavior of   the density for large $B.$ It asserts the existence of  
constants $a=a(L)$, $b=b(L)$ and $C=C(U,L, {\it v})$ $>0$, such that:
\begin{equation}\label{asymptotic}
N_{{\mathcal {L}}, {\it v}}(U, B)=C~B^{a}(\log B)^{b-1} (1+o(1)) \quad (B\rightarrow +\infty) 
\end{equation}
For nonsingular toric varieties, a suitable refinement of the conjecture was first
proved by Batyrev-Tschinkel \cite{batyrev}. Improvements in the error
term were then given by Salberger \cite{salberger}, and a bit later by
de la Breteche \cite{bretechetorique}.\par
Each of these works used a {\it particular} height function
$H_{\infty}$. The most important feature is its choice of $v_\infty\,,$  
which was defined as follows:
 $$\forall \y=(y_1,\dots, y_{n+1}) \in \Q_\infty^{n+1} (=\R^{n+1}), \qquad
v_\infty (\y) = \max_i |y_i|.$$
Although this is quite convenient to use for torics, there is no
fundamental reason why one should a priori limit the effort to prove the conjecture to this 
particular height function. 
In addition, since the existence of a precise asymptotic 
for one height function need not imply anything equally precise
for some other height function, a proof of  Manin's conjecture for
toric varieties and  
heights other than $H_\infty$ is not a consequence of the work cited above.\par 
In this paper we are interested in proving Manin's conjecture on
torics when a suitable homogeneous polynomial $P(X_1,\dots, X_{n+1})$
is used to define the metric's component $v_\infty$ at the archimedean
place. In this way,   a large number of  height functions $H_P$  on
$\P^n (\Q)$ are created. Our goal is then 
to compute the asymptotic density with respect to each $H_P$ of
any toric embedded in $\P^n(\Q)$.\par 
Since our problem focusses upon the choice of a polynomial
$P(X_1,\dots, X_{n+1})$ from the beginning, the dimension of the
ambient projective will equal $n$. We therefore start with  a toric 
embedded inside $\P^n(\Q) $ by a set of monomial defining equations as follows:
\begin{equation*}
X(\Ab) := \{(x_1:\dots:x_{n+1}) \in \P^{n}(\Q) \mid    
\prod_{j=1}^{n+1} \, x_j^{a_{i,j}} = 1 \, \forall i=1,\dots,l\}
\end{equation*} 
where $\bold A$ is an $l\times (n+1)$ matrix with entries in $\Z$,
whose rows   
$\a_i=(a_{i,1},\dots,a_{i,n+1})$ each satisfy the property that 
$\sum_{j=1}^{n+1} a_{i,j} = 0$. \\
The set of polynomials of interest in this paper consists of {\it
  generalized} polynomials $P=P (X_1, \dots, X_{n+1})$
  (i.e. exponents of monomials can be arbitrary nonnegative real numbers) with
positive  coefficients, elliptic on $[0, \infty)^{n+1}$ \footnote{$P$
  is elliptic on $[0, \infty)^{n+1}$ if its restriction to this set
  vanishes only at $(0,\dots, 0)$. In the following, the term
  ``elliptic" means ``elliptic on $[0, \infty)^{n+1}$\ ".} and homogeneous of degree 
$d>0.$ To each $P,$ we then consider the family of metrics $v_p$  on $\Q_p^{n+1}$ by setting,   $v_p(\y) = \max |y_j|_p$ for any $p < \infty$ and $\y = (y_1,\dots, y_{n+1})$ and   
\footnote{If $P$ is a quadratic form or
  is of the form $P=X_1^d+\dots+X_{n+1}^d$, then clearly the triangle
  inegality is verified. However, in order to preserve height under
  embedding, the definitions of metric used to build  heights 
(see for example (\cite{peyrebook}, chap. 2, \S 2.2) or \cite{peyre})
do not assume in general that the triangle inequality must hold.} 
 $$v_\infty (\y) = {\left(P(|y_1|,\dots,|y_{n+1}|)\right)}^{1/d} \qquad \forall \y=(y_1,\dots, y_{n+1}) \in \Q_\infty^{n+1}.$$
 A simple exercise
involving heights then shows that the height $H_P$ associated\footnote{Precisely, this is the
  height associated to the pair  $(\iota^*{\cal O}(1),  {\it v})$
  where   $\iota : X(\bold A)\hookrightarrow 
\P^n(\Q)$ is the canonical embedding,  $ {\cal O}(1)$  is the standard 
line bundle on $\P^n(\Q),$ and  ${\it v}=(v_p)_{p\leq \infty}$ is the
family of metrics. For simplicity this height is denoted $H_P.$}  to this family
of metrics is given as follows. For any $\x = (x_1:\dots:x_{n+1}) \in \P^{n}(\Q)$ such that  $gcd (x_i) = 1:$    
$$  H_P(\x):= {\left(P(|x_1|,\dots,|x_{n+1}|)\right)}^{1/d}\,.$$   

Very little appears to be known about the asymptotic density of
projective varieties for such heights $H_P$.   
A few earlier works are  \cite{peyre}, \cite{essouabrismf} or
\cite{swd}, but these are limited to very special choices of $P.$   
In particular,  
the reader should appreciate the fact that none of the general methods
developed to study the asymptotic 
density with respect to  $H_\infty$ can  be expected to    apply   to
any other height $H_P$.  
To be convinced of this basic fact, it is enough to look at the simple case where $X = \P^n(\Q)$ and  
$U=\{(x_1,\dots,x_{n+1}) \in  \P^n(\Q)\mid x_1 \dots x_{n+1} \neq 0\}$. 
In this case, an easy calculation  shows that 
\begin{eqnarray*}
N_{H_\infty}(U; t)&=&2^{n}\# \{\m \in \N^{n+1} \mid \max m_i
\leq t {\mbox { and }} gcd(m_i))=1\}\\
&\sim & \frac{2^n}{\zeta(n+1)}\# \{\m \in \N^{n+1} \mid \max m_i
\leq t \} \sim \frac{t^{n+1}}{\zeta(n+1)}\, . \quad (*1)
\end{eqnarray*}
On the other hand, if $P=P(X_1,\dots,X_{n+1})$ is  a homogeneous
elliptic polynomial of degree $d>0$ (for example $P=X_1^d+ \dots+X_{n+1}^d$), 
we have  
\begin{eqnarray*}
N_{H_P}(U; t)&=&2^{n}\# \{\m \in \N^{n+1} \mid P(\m) \leq t^d
{\mbox { and }} gcd(m_i)=1\}\\
&\sim& \frac{2^n}{\zeta(n+1)} \# \{\m \in \N^{n+1} \mid P(\m) \leq t^d \}\, .
\quad (*2) 
\end{eqnarray*}

The asymptotic  for the second factor  in 
$(*2)$ is given by 
$$\# \{\m \in \N^{n+1} \mid P(\m) \leq t^d \} \sim C(P;n) t^{n+1}, \quad  {\mbox { where }} \quad  
C(P;n)= \frac{1}{n+1} \int_{\Sp^{n} \cap \R_+^{n+1}} P_d^{-(n+1)/d}(\v)
d\sigma(\v),$$   $\Sp^ {n}$ denotes the unit sphere of $\R^n,$ 
and $d\sigma$ its Lebesgue measure.
This is a classical result of Mahler \cite{mahler}. It should be
evident to the reader that this  cannot be determined  
from the asymptotic   for the counting function of $(*1)$.\par 
Our  results are formulated  in terms of a  polyhedron in $[0,
\infty)^{n+1}\,,$ which can be associated   to 
$X(\bold A)$ in a natural way, and, in addition, the   maximal torus \\ 
$\ds U(\Ab )
:= \{(x_1:\dots:x_{n+1}) \in X(\Ab ) : x_1\dots x_{n+1} \neq 0
\}.$  

The first result, Theorem \ref{thMC} (see \S 3),  shows that if the diagonal
intersects the polyhedron in a compact face, then for any  height $H_P,$ as above,  there exist  
constants  $a=a(\bold A), b=b(\bold A)$ and $C=C({\bold A}, H_P)$ such that  
$$N_{H_P}(U(\bold A); t):=\# \{\x \in U(\bold A): H_P(\bold x) \leq t\}= C ~t^{a}
{\log(t)}^{b-1}\left(1+O\left((\log t)^{-1}\right))\right) {\mbox { as
  }} t \rightarrow \infty.$$ 
The constants $a$ and $b$  
 are also characterized quite simply  in terms of this
 polyhedron. The second part of Theorem \ref{thMC} refines this conclusion by
asserting that if the dimension of this face {\it equals} the 
dimension of $X(\bold A),$ then  $C > 0$. In this
event, we are also able to give an explicit expression for $C,$ the
form of which is a reasonable generalization of that given in
(*2).\par  
In order to derive Theorem \ref{thMC}, we establish in Theorem
\ref{thMCzeta} all the needed   analytical properties of the height zeta functions 
$\ds Z_{H_P}(U(\bold A); s):=\sum_{\x \in U(\bold A)} H_P(\x)^{-s}$.\par
Our second main result uses the fact that  we are able to  give a very precise description of the polyhedron  
  for the   class of hypersurfaces  $ \{ \x \in \P^{n}(\Q) \mid x_1^{a_1}\dots x_n^{a_n}
=x_{n+1}^{|\a|} \}.$  The particular case of the singular hypersurface 
$\ds \{\x \in \P^{n}(\Q) \mid x_1\dots x_n
=x_{n+1}^{n} \},$ provided with the height $H_\infty,$ was studied in several
papers (\cite{batyrev}, \cite{bretecheasterisque}, \cite{fouvry}, \cite{morozcubique},
\cite{eulerprod}). However, besides a result of
Swinnerton-Dyer \cite{swd} for  the singular cubic 
$X_1X_2X_3=X_4^3$ with a single height function {\it not} $H_\infty,$ nothing
comparable to our Theorem 2 appears to exist  in the  
literature (see Remark 4 in  \S 4.2). \par
The proof of these results are based on the new method I introduced in
\cite{essouamixte}. This allows one  to 
study analytic properties for 
``mixed zeta functions''  $Z(f;P;s)$ (see \S 4.1), and in particular,  to determine explicitly the  
principal part at its largest  pole (which is very important for proving the
asymptotic). Such zeta functions  combine together in one function
the multiplicative features of the variety and the additive nature of
the polynomial that defines the height $H_P$. Section \S 4.1 will be devoted to the adaptation of the general
results of \cite{essouamixte} 
to the case when $f$ is a multiplicative ``uniform'' function. Section \S 4.2 contains the proofs of our two theorems.

\section{Notations and preliminaries}
\subsection{Notations}
We gather in what follows some notations which will be used in this paper.
\be
\item 
$\N=\{1,2,\dots \}$,
$\N_0=\N \cup \{0\}$ and $p$ always 
denotes a prime number;
\item
The expression: $ f(\lambda,{\bf y},{\bf x}){\ll}_{{}_{{\bf y}}} g({\bf x})$
uniformly in ${ \bf x}\in X $ and ${\lambda}\in \Lambda$
means there exists $A=A({\bf y})>0$,
which depends neither on ${\bf x}$ nor $\lambda$, but could 
eventually depends on the parameter vector ${\bf y}$,
such that,   
$\forall {\bf x}\in X {\mbox { and }}\forall {\lambda}\in
{\Lambda}\quad |f(\lambda,{\bf y},{\bf x})|\leq Ag({\bf x}) $.
\item $f\asymp g$ means $f\ll g$ and $g\ll f$;
\item For any ${\bf x}=(x_1,..,x_n) \in \R^n$, we set   
$\|{\bf x}\|=\sqrt{x_1^2+..+x_n^2}$ and $|{\bf
  x}|=|x_1|+..+|x_n|$. We denote the canonical basis of $\R^n$ by 
$(\e_1,\dots,\e_n)$. The standard inert product on $\R^n$ is 
denoted by $\la .,. \ra$. We set also  
$\zerob=(0,\dots,0)$ and $\unb =(1,\dots,1)$.
\item 
We denote a vector in $\C^n$ $\s=(s_1,\dots,s_n)$, and write 
 $\s={\sigmab}+i{\taub},$
 where ${\sigmab } = (\sigma_1,\dots,\sigma_n)$ and  
 ${\taub }=(\tau_1,\dots,\tau_n)$ are the real resp. imaginary
 components of $\s$ (i.e. $\sigma_i=\Re(s_i)$ and 
$\tau_i=\Im(s_i)$ for all $i$). We also write 
$\la \x, \s\ra$ for $\sum_i x_i s_i$  if $\x \in \R^n, \s \in \C^n.$
\item Given $\alphab \in \N_0^n,$ we write $\X^{\alphab}$ for the
  monomial $X_1^{\alpha_1} \cdots X_n^{\alpha_n}$. 
For an analytic function  $h(\X)=\sum_{\alphab} a_{\alphab} \X^{\alphab}$, the set 
$supp(h):=\{\alphab \mid a_{\alphab} \neq 0\}$ is called the support of
$h$;
\item $f: \N^{n} \rightarrow \C$ is said to be 
  multiplicative if for all $m_1,\dots,m_n\in \N$ and 
  $m_1',\dots,m_n'\in \N$ satisfying 
$gcd\left(lcm\left(m_i\right),lcm\left(m_i'\right)\right)=1$ we have:
\\ 
$f\left(m_1 m_1',\dots , m_n m_n'\right)=f\left(m_1,\dots , m_n
\right) \cdot f\left(m_1',\dots ,m_n'\right)$;
\item Let $P(\X)= P(X_1,\dots, X_n)=\sum_{\alphab \in supp(P)}
  a_{\alphab} \X^{\alphab}$ be  
a generalized polynomial of degree $d$; i.e. $supp(P)$ is a finite subset of
$\R_+^n$ (not necessarly of $\N_0^n$). $P$ is said to be elliptic
if its homogenuous part of greater degree $P_d$ verifies:
$\forall \x \in \R^n_+ \setminus \{\zerob\}$, $P_d(\x)>0$;
\item Let $F$ be a meromorphic function on a domain ${\cal D}$ of
  $\C^n$ and let ${\cal S}$ be the support of its polar divisor. 
$F$ is said to be of moderate growth if there exist  
$a,b>0$ such that $\forall \delta >0$,  
 $F(\s) \ll_{\sigmab,\delta} 1+|\tau|^{a|\sigmab|+b}$ 
uniformly in $s=\sigmab+i\taub\in {\cal D}$ 
verifying $d(\s, {\cal S})\geq \delta $.
\ee
\subsection{Some preliminaries on the polyhedrons}
For the reader's convenience, notations and preliminaries about Newton
polyhedron that will be used throughout the article are assembled
here. For more details on the convex polyhedrons the reader can
consult for example the book \cite{roc}.\par
To any subset  $J$   of $[0, \infty)^n \setminus \{\zerob\},$ the boundary   
${\cal E}(J) = \partial \left(convex \ hull \big(\bigcup_{\alphab \in J}
  \alphab + \R^n\big)\right)$ is the Newton polyhedron generated by
$J$ The $J$ will be dropped from the notation if no confusion can 
occur as a result.
\be 
\item For any  $\a \in \R_+^n \setminus \{\zerob\}$  define
    $m(\a):=\inf_{\x \in {\cal E}} \la \a ,\x\ra$.  
Then $\a$ is a {\it polar vector} for the   face of ${\cal E} $ defined by
setting  
${\cal F}_{{\cal E}} (\a):=\{\x \in {\cal E} \mid \la \a , \x\ra =m(\a)\}.$
\item The faces  of ${\cal E}$ are: ${\cal F}_{{\cal E}} (\a)$ 
$(\a \in \R_+^n \setminus \{\zerob\})$.
We define a facet of $ {\cal E}$ as a face of maximal dimension.
\item
Let $F$ be a face of ${\cal E}$.
\be
\item 
The set $pol(F):=\{\a \in \R_+^n \setminus \{\zerob\} \mid
F={\cal F}_{{\cal E}} (\a)\}$ is called the polar set of $F.$ Its
elements are called polar vectors of $F$;
\item 
The set $pol_0(F):=\{\frac{\a}{m(\a)} \mid \a \in pol(F) 
{\mbox { and }} m(\a)\neq 0\}= \{\a \in pol(F) \mid m(\a)=1\}$ is
called the normalized polar set of the face $F.$ Its elements are
called normalized polar vectors of $F$;
\item If the face $F$ is not contained in a coordinate  hyperplane 
 then $\forall \a \in pol(F)$ $m(\a)\neq 0$. Moreover,   in this
 case we have $pol_0(F)  \neq \emptyset$. 
\item We set   ${\cal E}^0:= \partial \left( \bigcap_{\y \in {\cal E}} \{\x \in \R^n \mid \la \x ,\y \ra
  \geq 1 \}\right).$ This    convex set   is called the dual of ${\cal E}$; 
\ee
\item Let $A=\{\alphab^1,\dots,\alphab^q\}$ be a finite subset of $\R^n$. The convex cone of $A$ is $con (A):=\{\sum_{i=1}^q \lambda_i \alphab^i
\mid (\lambda_1,\dots,\lambda_q)\in \R_+^q \}$ and its  
(relative) interior is  
$con^* (A):=\{\sum_{i=1}^q \lambda_i \alphab^i
\mid (\lambda_1,\dots,\lambda_q)\in \R_+^{*q}\}$.
\ee 
\subsection{Construction of three important volume constants}
\subsubsection{First constant: The Sargos constant (\cite{sargosthese}, chap 3, \S 1.3)}
Let $P(\X)= \sum_{\alphab \in supp(P)} a_{\alphab} \X^{\alphab}$ be 
a generalized polynomial. We suppose that $P$ has 
positive coefficients and that it depends on all the variables $X_1,\dots,X_n$.
We denote by 
${\cal E}^{\infty}(P):=\partial \left(conv(supp(P))-\R_+^n\right)$ its
Newton polyhedron at infinity. 
Let $G_0$ be the smallest face of ${\cal E}^{\infty}(P)$ which meets
the diagonal $\Delta =\R_+ \unb $.
We denote by $\sigma_0=\sigma_0(P)$ the unique positive real number
$t$ which verifies $t^{-1}\unb \in G_0$ and we set $\rho_0 =\rho_0 (P):=codim G_0$. 
By coordinates permutation one can suppose that $\oplus_{i=1}^{\rho_0}\R \e_i
\oplus \overrightarrow{G_0}=\R^n$ and that 
$\{\e_i \mid G_0=G_0-\R_+ \e_i\}=\{e_{m+1},\dots,\e_n\}$.\\
Let $\lambdab_1,\dots,\lambdab_N$ be the normalized polar vectors of the facets
of ${\cal E}^{\infty}(P)$ which meet $\Delta$.\\
Set $P_{G_0}(X)=\sum_{\alphab \in G_0} a_{\alphab} \X^{\alphab}$ and 
$\Lambda =Conv\{\zerob,\lambdab_1,\dots,\lambdab_N, \e_{\rho_0+1},\dots,
\e_n\}$.\\
The {\bf Sargos constant} associated to $P$ is: 
$$A_0(P):= n!~Vol(\Lambda)~\int_{[1,+\infty[^{n-m}}\left(\int_{\R_+^{n-\rho_0}}P_{G_0}^{-\sigma_0}
(\unb,\x, \y) ~d\x\right) d\y.$$
\subsubsection{Second constant: The volume constant}
Let $I$ be a finite subset of $\R_+^r \setminus \{\zerob\}$, 
$\u=(u(\betab)_{\betab \in I})$ a finite sequence of positive integers
and $\b =(b_1,\dots,b_r)\in \R_+^{*r}$.
To these three data, we associate the generalized polynomial  
$P_{(I;\u;\b)}$ with $q:=\sum_{\betab \in I} u(\betab)$ variables, 
in the following way:\\
We define $\alphab^1,\dots,\alphab^q$ by:
$\{\alphab^i \mid i=1,\dots,q\}=I$ and
$\forall \betab \in I$ 
$\#\{i\in \{1,..,q\} \mid \alphab^i=\betab \}=u(\betab)$ (i.e. the
family $(\alphab^i)$ is obtained by repeating $u(\betab)$ times each element $\betab$ of $I$).\\
We define the vectors $\gammab^1,\dots,\gammab^r$ 
of $\R_+^q$ by: $\gamma_i^k=\alpha_k^i$ \quad $\forall i=1,\dots,q$ and $\forall k=1,\dots,r.$ \\
Finally, we set  $P_{(I;\u;\b)}(\X):=\sum_{i=1}^r b_i
\X^{\gammab^i}$.\\
The {\bf volume constant} $A_0(I;\u;\b)$ associated to $I$, $\u$, $\b$ is the Sargos
constant associated to the polynomial $P_{(I;\u;\b)}$; i.e. 
$$A_0(I;\u;\b):=A_0\left(P_{(I;\u;\b)}\right)>0.$$
\subsubsection{The third constant: The mixed volume constant}
Let  
$P=P(X_1,\dots, X_n) =b_1 \X^{\gammab^1}+\dots+b_r
\X^{\gammab^r}$ be a generalized polynomial with positive coefficients. We set $\b =(b_1,\dots,b_r)\in
\R_+^{*r}$.\\ 
Let $\T=(I,\u),$ where $I$ is a finite subset of 
$\R_+^n\setminus \{0\},$ and $\u =\left(u(\betab)\right)_{\betab \in I}$ is 
a finite family of positive integers.\\
We associate to $\T$ and $P$ the following objects:
\be
\item $n$ elements $\alphab^1,\dots,\alphab^n$ of
  $\R_+^r$ defined by:
$\alpha^i_j =\gamma^j_i$ $\forall i=1,\dots,n$ and  $\forall
j=1,\dots,r$;
\item $\mu(\T;P;\betab):=\sum_{i=1}^n \beta_i \alphab^i$ for all  
$\betab \in I$ \ \ and \ \ $I_{\T,P}=\{\mu(\T;P;\betab) \mid \betab \in
I \}$; 
\item $\u_{\T,P}=\left(u_{\T,P}(\etab)\right)_{\etab \in I_{\T,P}}$ \ \ where \ \  
$u_{\T,P}(\etab)=\sum_{\{\betab \in I;
  ~\mu(\T;P;\betab)=\etab\}} u(\betab)$  
$\forall \etab \in I_{\T;P}$;
\ee
We define finally the {\it mixed volume} constant associated to $\T$
and $P$ by:
$$
A_0(\T;P)=A_0(I_{\T,P};\u_{\T,P};\b)\quad {\mbox { where}}
$$
$A_0(I_{\T,P};\u_{f,P};\b)>0$ is the volume constant from \S 2.3.2 that is 
  associated to $I_{\T,P};\u_{\T,P}$ and $\b$.

\section{Statements of main results}
Let $\bold A$ a $l\times (n+1)$ matrix with entries in $\Z$,
whose rows  $\a_i=(a_{i,1},\dots,a_{i,n+1})$ each satisfy the property that 
$\sum_{j=1}^{n+1} a_{i,j} = 0$. 
We consider the projective toric  varieties defined by: 
\begin{equation}\label{toriquedef}
X(\Ab) := \{(x_1:\dots:x_{n+1}) \in \P^{n}(\Q) \mid    
\prod_{j=1}^{n+1} \, x_j^{a_{i,j}} = 1 \, \forall i=1,\dots,l\}
\end{equation} 
We assume, without loss of generality, that the rows
$\a_i$ $(i=1,\dots,l)$ are linearly independent over $\Q$. It follows that: 
\begin{equation}\label{dimensionXA}
rank(\Ab )=l {\mbox { and }} dim X(\Ab )= n -rank(\Ab ) = n -l.
\end{equation}
Denote by 
$\ds U(\Ab )
:= \{(x_1:\dots:x_{n+1}) \in X(\Ab ) : x_1\dots x_{n+1} \neq 0
\}$ the maximal torus of $X(\Ab )$.\\
Define also  
\begin{eqnarray}\label{ca} 
T(\bold A) 
&:=& \left\{\bold \nu \in \N_0^{n+1} \mid \bold A(\bold \nu) = \bold 0 {\mbox {
    and }}\prod_i \nu_i = 0\right\} {\mbox { and }} T^*(\bold A)= T(\bold
A)\setminus \{\zerob\};\nonumber \\
c(\bold A) &:=& \frac{1}{2} \#
\big\{ {\bf \epsilon } \in \{-1,+1\}^{n+1} \mid
\prod_{j=1}^{n+1}  \epsilon_j^{a_{i,j}}=1 \, \forall i \big\}
\end{eqnarray}
Define ${\cal E}(\bold A) = {\cal E}\left(T^*(\bold A)\right),$ and set 
${\cal F}_0(\bold A)$ to denote the smallest face of ${\cal E}(\bold
A)$ that meets the diagonal. We then introduce the following: 
\be
\item $\rho(\bold A):=\#\left( {\cal F}_0(\bold A)\cap T^*(\bold A) \right) -
dim \left({\cal F}_0(\bold A)\right)$;
\item ${\cal E}^0 (\bold A)  = \partial \left( \bigcap_{\nub \in {\cal E}(\bold A)} \ \{\x
  \in \R^n \ \mid \ \la \x, \nub \ra \geq 1 \} \right)$ ( the dual of  ${\cal E}(\bold A)$);
\item $\iota (\bold A):= \min \{|\cb| \ \mid \  
  \cb \in {\cal E}^0 (\bold A)\cap \R_+^n\}$.
\ee
Fix now a (generalized) polynomial $P=P(X_1,\dots,X_{n+1})$ with
positive coefficients and assume that $P$ is 
elliptic and homogeneous of degree $d>0$. Denote by $H_P$ the height 
of $\P^{n}(\Q)$ associated to $P$. We introduce also the following notations: 
\be
\item Writing $P$ as a sum of monomials $P(X)= b_1 X^{\gammab^1}+\dots+b_r
X^{\gammab^r}\,,$  we set 
$\b =(b_1,\dots,b_r)\in \R_+^{*r}$ 
 \item Defining  $\alphab^1,\dots,\alphab^{n+1} \in \R_+^r $ to be the
   row  vectors  of the  matrix that equals the transpose of the
   matrix with  rows     $\gammab^1,..., \gammab^r,$ we set 
$I^*(\bold A)=\{ \sum_{i=1}^n \beta_i \alphab^i \mid \betab \in  
{\cal F}_0(\bold A)\cap T^*(\bold A)\}.$   
\ee
We can now state the first result as follows.
\begin{Theoreme} \label{thMC}
If the Newton polyhedron ${\cal E}(\bold A)$
 has a compact face which meets the diagonal,  
then there exists a polynomial $Q$ of degree at most $\rho (\bold A)
-1$ and  $\theta >0$
such that as  $t\rightarrow +\infty $:
$$N_{H_P}(U(\bold A);t):=\#\{M \in U(\bold A)\mid  H_P(M)\leq t\} 
=t^{\iota (\bold A)} Q(\log(t)) + O\left(t^{\iota (\bold A)-\theta}\right).$$
If we assume in addition that 
$dim\left({\cal F}_0(\bold A)\right)=dim X(\bold A)=n -l$, 
then $Q$ is a nonzero polynomial of degree   
$\rho (\bold A) -1$ and:
$$N_{H_P}(U(\bold A);t) = C\left({\bold A}; H_P\right)~ t^{\iota (\bold A)} (\log t)^{\rho (\bold A) -1} ~
\left(1+O\left((\log t)^{-1}\right)\right)$$
where 
$$ 
C\left({\bold A}; H_P\right):= \frac{c(\bold A)~ d^{\rho (\bold A)}
    A_0(I^*(\bold A);\unb;\b)}{\iota  (\bold A) ~ (\rho  (\bold A)-1)!}~
\prod_p \big[\big(1-\frac{1}{p}\big)^{\# {\cal F}_0(\bold A)\cap 
T^*(\bold A)}
\big(\sum_{\nub \in T(\bold A)} \frac{1}{p^{\la \nub, \cb\ra }}\big)
\big] >0,$$  
 $c({\bf A})$ is defined by (\ref{ca}),  
$A_0(I^*(\bold A); \unb;\b)$ is the volume constant (see \S 2.3.2),   
and $\cb$ is any\footnote{The constant 
  $C\left({\bold A}; H_P\right)$ does not depend on this choice.} 
normalized polar vector of the face ${\cal F}_0(\bold A)$.
\end{Theoreme}

Theorem \ref{thMC} follows from the analytic properties of heights zeta
functions established in the following.
\begin{Theoreme} \label{thMCzeta}
If the Newton polyhedron ${\cal E}(\bold A)$
has a compact face which meets the diagonal,  
then the height zeta function 
$\ds \s \mapsto Z_{H_P}(U(\bold A); s):=\sum_{M \in U(\bold A)}
H_P^{-s}(M)$ is   holomorphic  
in the half-plane \\
$\{s\in \C \mid \sigma >\iota (\bold A)\},$ and there exists $
\eta >0$ such that $s \mapsto Z_{H_P}(U(\bold A); s)$ has meromorphic continuation
with moderate growth to the half-plane 
$\{\sigma > \iota (\bold A) - \eta\}$ with only one possible pole at 
$s=\iota(\bold A)$ of order at most $\rho (\bold A).$

If we assume in addition that   
$dim\left({\cal F}_0(\bold A)\right)=dim X(\bold A)=n -l$, 
then $s=\iota(\bold A)$ is indeed a pole of order $\rho  (\bold A)$ and 
$\ds Z_{H_P}(U(\bold A); s) \sim_{s\rightarrow \iota (\bold A)} 
\frac{C_0\left({\bold A}; H_P\right)}{\left(s-\iota(\bold A)\right)^{\rho
    (\bold A)}}, \quad {\mbox {where }}$
$$ 
C_0\left({\bold A}; H_P\right):= C\left({\bold A}; H_P\right)~\iota (\bold A) ~ (\rho  (\bold A)-1)! >0,$$  
($C\left({\bold A}; H_P\right)$ is the constant volume defined in Theorem
\ref{thMC} above). 
\end{Theoreme}

Theorem \ref{thMC} (or equivalently  Theorem \ref{thMCzeta}) is a general result that applies to any projective
toric variety. Our second result applies Theorem \ref{thMC} to a particular class of
toric hypersurfaces that correspond to a class of problems from multiplicative number theory.\par
Let $n\in \N$ $(n\geq 2)$ and $\a=(a_1,\dots,a_n)\in \N^{n}$. Set 
$q=|\a|=a_1+\dots +a_n$. \\
Consider the  hypersurface: 
$\ds
X_n(\a)= \{\x \in \P^n (\Q) \mid  
x_1^{a_1}\dots x_n^{a_n} = x_{n+1}^q \}$ 
with torus:\\
$\ds U_n(\a) = \{\x \in X_n(\a)\mid 
x_1\dots x_n \neq 0\}.$ 

Let $P =P(X_1,\dots ,X_{n+1})$ be a generalized polynomial 
  as in Theorem 1. Set\\
   ${\tilde P}(X_1,\dots,X_n):=P(X_1,\cdots,X_n,\prod_j X_j^{a_j/q}),$  and write 
${\tilde P}(\X)=\sum_{k=1}^r b_k \X^{\gammab^k}$. \\ 
Define:
\begin{eqnarray*}
L_n(\a) &:=& \left\{\r \in 
\N_0^n~~;~~   
q | \la \a, \r\ra {\mbox { and }} r_1\dots r_n =0 \right\}\setminus \{\zerob\};\\
{\cal E}(\a)&:=& {\cal E}(L_n(\a)) ;\\
{\cal F}_0(\a) &=&  {\mbox {the smallest face of }} {\cal E}(\a)
{\mbox { which meets
the diagonal }} \Delta =\R_+ \unb ;\\
J_n(\a) &:=& L_n(\a)\cap {\cal F}_0(\a) \quad {\mbox { and }} \quad \rho (\a)=\#
\left(J_n(\a)\right) -n+1;\\
J_n^*(\a)&:=&\{\sum_{i=1}^{n} \beta_i \alphab^i \mid \betab 
\in J_n(\a)\};\\
c(\a)&:=& \frac{1}{2}\#\left\{(\eps_1,\dots ,\eps_{n+1})\in \{-1,+1\}^{n+1}
\mid \eps_1^{a_1}\dots \eps_n^{a_n}=\eps_{n+1}^q \right\}.
\end{eqnarray*}
Applying Theorem 1 we obtain the following:
\begin{Theoreme}\label{CMa}
Let  $ \bold c  \in \R_+^{*n}$ be a normalized polar vector of the
face ${\cal F}_0(\a)$.  
There exists a polynomial $Q$ of degree at most 
$\rho (\a)-1$ and $\theta >0$ such that:
$$
N_{H_{P}}(U_n(\a);t)=\#\{M \in U_n(\a) \mid  H_{P}(M)\leq t\}
= t^{|\bold c|} Q(\log t) + O\left(t^{1-\theta}\right).$$
If we assume in addition that ${\cal F}_0(\a)$ is a facet of 
the polyhedron ${\cal E}(\a)$, then $Q\neq 0$, 
deg$Q = \rho (\a)-1$ and  
$$ N_{H_{P}}(U_n(\a); t)= C\left(\a ; H_P\right) ~t^{|\bold c|} {(\log t)}^{\rho (\a)-1} 
\times \bigg(1+O\big({(\log t)}^{-1}\big)\bigg)  $$
where:
$$ C\left(\a ; H_P\right):=
\frac{c(\a)~d^{\rho (\a)}~ A_0\left(J_n^*(\a);\unb;
    \b\right)}{|\bold c| \cdot (\rho (\a)-1)!}\, \cdot \, C_n (\a)$$
$$  C_n (\a):= \prod_p \left(
\bigg(1 - \frac{1}{p}\right)^{\rho(\a)+n-1}
\, \bigg(\sum_{\nub \in \N_0^n;~~q~|~\la \a ,\nub \ra \atop 
\nu_1 \dots \nu_n =0} p^{-\la \cb(\a) ,\nub \ra}\bigg)\bigg) >0\,,$$
and $A_0\left(J_n^*(\a);\unb;\b\right)>0$ is the volume
 constant.
\end{Theoreme}
{\bf Remark 1:} An interesting question is to determine the precise
set of  exponents $a_1,\dots, a_n \ge 1$ (for given $q$ and $n$) such
that ${\cal F}_0(\a)$ 
is a facet of ${\cal E}(\a).$ It seems   reasonable to  believe that
the complement of this set is thin in a suitable sense (when $q$ is allowed to be arbitrary).\par
{\bf Remark 2:} {\it If $\a=(a_1,\dots,a_n)\in \N^{n}$ satisfies the property that  each  
$a_i$ divides  $q = a_1+\dots+a_n$, then  
${\cal F}_0(\a) =conv\left(\{\frac{q}{a_i} \e_i \mid
  i=1,\dots,n\}\right)$. Thus,  ${\cal F}_0(\a)$ is a facet of 
${\cal E}(\a)$ and  the more precise second part of  Theorem \ref{CMa} applies.}\par 
Manin's conjecture 
had been proved for the height $H_\infty$ and the particular surface  $X_3 (\unb)$ in several earlier works (see 
\cite{fouvry}, \cite{bretecheasterisque}, \cite{salberger}). More recently,  the article \cite{eulerprod}
extended these earlier results to   $X_n(\unb)$ for any $n \ge 3$ (but
only used  $H_\infty$). Theorem \ref{CMa} should  therefore  be understood as a natural  generalization of  all   
these earlier results.\par 
{\bf Remark 3:} {\it When $\a = \unb,$ it is interesting to compare the
constant $C(\unb ; H_\infty)$ with $C (\unb ; H_{P_d})$
where   $P_d:=\sum_{k=1}^{n+1} X_k^{d}$. Theorem \ref{CMa} implies:   
$$ C\left( \unb ; H_{P_d}\right) := \frac{2^{n-1} ~d^{{2n-1\choose
      n}-n}~ A_0(P_{d,n})}{({2n-1\choose n}-n-1)!} 
\left[\prod_{p}\left( {\left(1 - \frac{1}{p}\right)}^{{2n-1\choose n}}
\, \left(\sum_{k=0}^\infty  \frac{{(k+1)n-1\choose n-1}}{p^k}\right) \right)\right]>0,$$
where  $P_{d,n}$ is defined as follows. First, note that  $ J_n(\unb)=\{ \r \in \N_0^n \mid r_1+\dots+r_n =n 
{\mbox { and }} r_1\dots r_n =0\}.$ Define    
$\rho_n:=\# J_n(\unb) ={2n-1\choose n} -1$. 
Let $\left(m_{i,j}\right)$ be the $(n+1) \times \rho_n$ 
matrix whose columns are the vectors 
$(d r_1,\dots,d r_n, d)$ $(\r \in J_n(\unb)).$ Now define  
$$P_{d,n}(X_1,\dots, X_{\rho_n}):=\sum_{i=1}^{n+1} X_1^{m_{i,1}}\dots
X_{\rho_n}^{m_{i,\rho_n}}.$$ 
The constant $A_0 (P_{d,n})$ is the Sargos constant (\S 2.3.2).\\   
As noted in the Introduction, it is not at all evident that the expression for  
$C\left( \unb ; H_\infty\right)$ (see  \cite{eulerprod}) determines
that for  $C\left( \unb ; H_{P_d} \right)$ given above.  For
example in the simplest case $n=2$, a direct calculation gives 
$C\left( \unb ; H_\infty \right) = \frac{12}{\pi^2}\,,$ 
whereas Theorem 2 gives:\\   
$C\left( \unb ; H_{P_d}\right)=\frac{6}{\pi^2} 
\int_0^{\frac{\pi}{2}} \left(\cos^{2d}(\theta)+\sin^{2d}(\theta)+
    \cos^{d}(\theta)\sin^{d}(\theta)\right)^{-1/d} ~d\theta$.}

\section{Proofs }  
The proofs of Theorems  \ref{thMC}, \ref{thMCzeta} and \ref{CMa} are simple consequences of the main result
proved in \S 4.1. This derives basic analytic properties of a ``mixed
zeta function'' that combines both the multiplicative and additive
features of our underlying counting problem into one generating
function.\par 
For simplicity, we will work in \S 4.1 (as well as the  Appendix) with
functions $f$ and polynomials $P$ of $n$ variables instead of $n+1$
variables as in the remainder of the article.  The reader should be alert to this shift in index when reading \S 4.2.
\subsection{Mixed zeta functions}  
Let $f:\N^{n} \mapsto \C$ be an arithmetical function and let $P$ 
be a (generalized) polynomial of degree $d$. Consider the ``mixed zeta function"
defined formally for $s\in \C$ by 
$$  Z(f;P;s):= \sum_{\m \in \N^n}\frac{f(m_1,\dots,m_n)}{P^{s/d}(m_1,\dots,m_n)}\,.$$ 
In \cite{essouamixte}  I proved meromorphic continuation and
determined explicitly the principal part at the first pole, assuming   that  certain reasonable 
properties were satisfied by   $f$ and $P.$  For the reader's convenience,  a summary
of this is given  in the Appendix.\par 
In this section we will assume $f$ is a  uniform multiplicative 
function   and $P$ is elliptic and homogeneous. These conditions
suffice to prove our theorems.     
\begin{Definition}\label{fmh}
A multiplicative function $f:\N^{n} \mapsto \N_0$ is said to be
uniform if there exists a function $g_f:\N_0^{n} \mapsto \N_0$
and two constants $M,C>0$ such that for all prime numbers  $p$ and  
all  $\nub \in \N_0^n$,  $f(p^{\nu_1},\dots,p^{\nu_n})=g_f(\nub)\leq C \left(1+|\nub|\right)^M$. 
\end{Definition}
We fix a uniform multiplicative  
function $f:\N^{n} \mapsto \N_0$  throughout the rest of \S 4.1 and
write $g$ in place of $g_f$.    
We set  $S^*(g)=\{\nub \in \N_0^n\setminus
\{\zerob\} \mid g(\nub)\neq 0\}$ and assume that $S^*(g)\neq \emptyset$.
We then define:
\be
\item
${\cal E}(f):={\cal E}(S^*(g)) =$   the Newton  polyhedron determined by   $S^*(g)$;
\item ${\cal E}(f)^o:=$  the dual of  ${\cal E}(f)$;
\item $\iota (f):=\min \left\{|\cb| \mid
    \cb \in {\cal E}(f)^o \cap
    \R_+^n \right\}$ (  the ``index" of $f$); 
\item ${\cal F}_0 ( f ) := $ the smallest face of ${\cal E}(f)$ which meets the diagonal. We denote its set of 
normalized polar vectors by  
$pol_0\left({\cal F}_0 (f) \right);$  
\item $\T:=(I;\u)$ where $I:={\cal F}_0 (f) \cap S^*(g)$
and $\u= \left(g(\betab)\right)_{\betab \in I}$.
\ee
For $f, P$ as above, we will deduce our main results from the
following Theorem. 
\begin{Theoreme}\label{ftfpourapplication}
If the face ${\cal F}_0(f)$ is compact, then 
$\ds Z(f;P;s):=\sum_{\m \in \N^{n}}
\frac{f(m_1,\dots,m_n)}{P(m_1,\dots,m_n)^{s/d}}$ is {\it holomorphic}
in the half-plane $\{s: \sigma >\iota(f)\}$ and there exists $
\eta >0$ such that $s \mapsto Z(f;P;s)$ has meromorphic continuation
with moderate growth to the half-plane 
$\{\sigma > \iota (f) - \eta\}$ with only one possible pole at 
$s=\iota(f)$ of order at most $\rho_0 (\T):=\sum_{\betab \in I}
g(\betab) -rank(I)+1$.\par
If, in addition, we assume that   
 $dim \,{\cal F}_0(f) = rank(S^*(g))-1$, then:  
$s=\iota(f)$ is indeed a pole of order $\rho_0 (\T)$ and 
$$Z(f;P;s)\sim_{s\rightarrow \iota(f)} \frac{C_{ari}(f)~ d^{\rho_0 (\T)}~
A_0(\T;P)}{\left(s-\iota(f)\right)^{\rho_0 (\T)}}, \quad {\mbox {where }}$$
$\ds C_{ari}(f)= \prod_p\left(
\left(1-\frac{1}{p}\right)^{\sum_{\nub \in I} g(\nub)}~
\left(\sum_{\nub \in \N_0^n} \frac{g(\nub)}{p^{\la \nub,
      \cb\ra}}\right)\right)>0$ 
($\cb$ is any\footnote{The constant 
  $C_{ari}(f)$ doesn't depend on the choice of $\cb$} normalized polar vector of  
${\cal F}_0 (f)$  \\
and $A_0(\T;P)>0$ is the {\it mixed volume} constant
defined in \S 2.3.
\end{Theoreme}

{\bf Remark:}
The assumption $dim \, {\cal F}_0 (f) = rank(S^*(g))-1$ 
is automatically satisfied  if for example the face 
${\cal F}_0 (f)$ is a facet of ${\cal E}(f)$. 
\par
By a simple adaptation of a standard tauberian argument of Landau (see
for example \cite{essouabrismf}, Prop. 3.1)), 
we obtain the following arithmetical consequence:
\begin{Corollaire}\label{applicationmulti}
If the face ${\cal F}_0 (f)$ is compact, then there exist 
$\delta >0$ and a polynomial $Q$ of degree at most  
$\rho (f;P)-1$ such that:
$$N(f;P;t):=\sum_{\{\m \in \N^{n}; P(\m)^{1/d} \leq t\}}
f(m_1,\dots,m_n) = t^{\iota(f)} Q(\log t) + O(t^{\iota(f)-\delta})
{\mbox { as }}  t\rightarrow
+\infty.$$
If we assume in addition that 
$dim \,{\cal F}_0 (f)=rank(S^*(g))-1$, then $Q\neq 0$,  
 deg$Q=\rho (f;P)-1$ and 
$$N(f;P;t):= \frac{C_{ari}(f) d^{\rho_0 (\T)}
A_0(\T;P)}{\iota(f) \cdot  \left(\rho_0(\T)-1\right)!}~t^{\iota(f)} (\log t)^{\rho_0(\T)-1} ~\left(1+O\left((\log
  t)^{-1}\right)\right) {\mbox { as }}  t\rightarrow
+\infty.$$
\end{Corollaire}

\subsubsection{Proof of Theorem \ref{ftfpourapplication}}
Using the definitions introduced in \S 2.2, I will first give a lemma from convex analysis.
\begin{Lemme}\label{lienfaceetra}
Let $I$ be a nonempty subset of $\R_+^n \setminus \{\zerob\}$. Set ${\cal
  E}(I) $ to be its Newton polyhedron and denotes by ${\cal E}^o(I)$
its dual (see \S 2.2). Let $F$ be a face of  ${\cal E}(I)$ that is not contained
in a coordinate  hyperplane  and $\cb\in pol_0(F)$ a
normalized polar vector of $F$. Then: 
$F$ meets the diagonal if and only if 
$|\cb|=\iota (I)$ where  
$\iota (I) =  \min\{|\alphab |;~
\alphab \in {\cal E}^o(I) \cap \R_+^n\}$.
\end{Lemme}
{\bf Proof of Lemma \ref{lienfaceetra} }\\
We first note that the definition of the normalized polar vector
implies that $\cb \in {\cal E}(I)^o \cap \R_+^n$.\\
$\bullet$ Assume first that the diagonal $\Delta$ meets the face $F$ 
of the Newton Polyhedron ${\cal E}(I)$.
Therefore, there exists $t_0 >0$ such that $\Delta \cap F =\{t_0 \unb\}$.
Let $\alphab^1,\dots,\alphab^r \in I \cap
F$ and let $J$ a subset (possibly empty) of $\{1,\dots,n\}$
such that $F = convex \ hull \{\alphab^1,\dots,\alphab^r\}+con\{\e_i \mid i\in J\}$.
Thus there exist $\lambda_1,\dots,\lambda_r \in \R_+$
verifying $\lambda_1+\dots +\lambda_r =1$ and a finite family $(\mu_i)_{i\in J}$
of elements of $\R_+$ such that:
\begin{equation} \label{t0un}
t_0 \unb = \sum_{i=1}^r \lambda_i \alphab^i + \sum_{j\in J} \mu_j
\e_j.
\end{equation} 
But $\cb$ is orthogonal to the vectors $\e_j$ $(j\in J)$ and   
$\la \cb , \alphab^i \ra =1$ for all $i=1,\dots,r$. Thus it follows
from the relation
(\ref{t0un}) that 
$t_0 |\cb|=\la \cb , t_0 \unb\ra =
\sum_{i=1}^r \lambda_i \la \alphab^i , \cb \ra + \sum_{j\in J}
\mu_j  \la \e_j , \cb \ra = \sum_{i=1}^r \lambda_i =1$. Consequently  
$|\cb|=t_0^{-1}$. \par
Relation (\ref{t0un}) implies also that for all 
$\b \in {\cal E}(I)^o \cap \R_+^n$: 
$$|\b|=t_0^{-1}~ \la \b , t_0 \unb\ra =
t_0^{-1}~ \sum_{i=1}^r \lambda_i \la \alphab^i , \b \ra + \sum_{j\in J}
\mu_j b_j 
\geq t_0^{-1}~ \sum_{i=1}^r \lambda_i \la \alphab^i , \b \ra \geq
t_0^{-1}~ \sum_{i=1}^r \lambda_i =t_0^{-1}=|\cb|.$$
This implies $|\cb|=\iota (I).$  \par
$\bullet$ Conversely assume that $|\cb|=\iota (I)$. 
We will show that $\Delta \cap F \neq \emptyset\,.$\\ 
Let $G$ be a face of ${\cal E}(I)$ which meets the diagonal 
$\Delta$. Since $G$ is not included in the coordinate hyperplanes, it
has a normalized polar vector $\a \in pol_0(G)$.  
Moreover there exist $\betab^1,\dots,\betab^k \in I \cap G$ and 
$T$ a subset (possibly empty) of $\{1,\dots,n\}$ such that\\ 
$G  = convex \ hull \{\betab^1,\dots,\betab^k\}+con\{\e_i \mid i\in T\}$.\\
The proof of the first part shows then that $|\a|=\iota (I)$ and that
there exist $\nu_1,\dots,\nu_k \in \R_+$
verifying $\nu_1+\dots +\nu_k =1$ and a finite family $(\delta_j)_{j\in T}$
of elements of $\R_+$ such that if we set 
$t_0 :=|\a|^{-1}$,then:
\begin{equation} \label{t0unbis}
t_0 \unb = \sum_{i=1}^k \nu_i \betab^i + \sum_{j\in T} \delta_j
\e_j \in G\cap \Delta.
\end{equation} 
Set $T':=\{ j\in T \mid \delta_j \neq 0\}$. 
It follows from relation (\ref{t0unbis}) that 
$$|\cb|=t_0^{-1}~ \la \cb , t_0 \unb\ra =
t_0^{-1}~ \sum_{i=1}^k \nu_i \la \betab^i , \cb \ra + \sum_{j\in J'}
\delta_j \la \e_j, \cb\ra  
\geq t_0^{-1}~ \sum_{i=1}^k \nu_i \la \betab^i , \cb \ra \geq
t_0^{-1}~ \sum_{i=1}^k \nu_i =t_0^{-1}=|\a|.$$
But $|\a|=\iota(I)=|\cb|$ so the intermediate inequalities must be equalities. This clearly forces 
$\la \betab^i,\cb \ra =1$ $\forall i=1,\dots,k$ and $\la \cb , \e_j\ra =0$
$\forall j\in T'$.
Relation (\ref{t0unbis}) implies then that:\\ 
$\la t_0 \unb, \cb \ra =
\sum_{i=1}^k \nu_i \la \betab^i,\cb\ra + \sum_{j\in J'} \delta_j
\la \e_j,\cb \ra = \sum_{i=1}^k \nu_i =1$.\\
Since $t_0 \unb \in {\cal E}(I)$,  it follows from the preceding discussion  that   
$t_0 \unb \in F$ and therefore $\Delta \cap F \neq \emptyset$. 
This finishes the proof of Lemma \ref{lienfaceetra}.\cqfd \par

The second step needed to prove Theorem \ref{ftfpourapplication} 
is a proposition that shows that the assumptions, introduced 
in the statement of Theorem A (see Appendix), are in fact satisfied 
whenever $f$ is a uniform multiplicative function. Using the notations
introduced after Definition 1, this proposition is as follows. 

\begin{Proposition}\label{propofg}
Let ${\cal B}(g)$ be the set of   $\betab \in S^*(g)$ which are 
in at least one compact face of the polyhedron ${\cal E}(f)$. 
Then: 
$\ds {\cal M}(f;\s)=\sum_{m_1,\dots,m_n \geq 1}
\frac{f(m_1,\dots,m_n)}{m_1^{s_1}\dots m_n^{s_n}}$ converges
absolutely in 
$${\cal D}_f:=\{\s \in \C^n \mid \Re \la
\s, \betab \ra > 1  ~\forall \betab \in
{\cal B}(g)\}\cap \{ \s \in \C^n \mid \Re s_i > 0 ~\forall
i=1,\dots,n\}.$$
Moreover for any $\cb \in pol_0 ({\cal F}_0 (f))$, ${\cal M}(f;\s)$ converges absolutely  
in $\{\s\mid \Re(s_i) > c_i ~\forall i\}$ and there exists $\eps_0>0$
such that  
$$\s\mapsto H_{\cb}(f;\s):= \left(\prod_{\betab \in {\cal F}_0 (f) \cap S^*(g)} 
(\la \betab, \s \ra)^{g(\betab)}\right)~{\cal M}(f;\cb+\s)$$ has a
holomorphic continuation with moderate growth to $\{\s \in \C^n \mid
\forall i ~\Re(s_i) >-\eps_0\}$ and satisfies    
$$H_{\cb}(f;\zerob) =\prod_p
\left(1-\frac{1}{p}\right)^{\sum_{\nub \in {\cal
      F}_0 (f) \cap S^*(g) }g(\nub)}~\left(\sum_{\nub \in
    \N_0^n} \frac{g(\nub)}{p^{\la \nub,\cb\ra}}\right) >0.$$
\end{Proposition}

{\bf Proof of Proposition \ref{propofg}:}\\
For all  $\eps \in \R$, set 
$U_\eps :=\{\s \in \C^n \mid \Re(s_i) >\eps ~\forall i=1,\dots,n\}$.\\
$\bullet$ Let $\eps >0$. Fix 
$M:=\left[\frac{8}{\eps}\right]+1\in \N$. 
It is easy to see that we have uniformly in $p$  and   $\s
=\sigmab+i\taub \in U_\eps $:
$$
\sum_{|\nub| \geq M+1}\frac{g(\nub)}{p^{\la \sigmab ,\nub\ra}}
\ll_\eps \sum_{|\nub| \geq M+1}\frac{|\nub|^D}{p^{\eps |\nub|}}
\ll_\eps \frac{1}{p^{\eps (M+1)/2}}
\sum_{|\nub| \geq M+1}\frac{|\nub|^D}{2^{\eps |\nub|/2}}
\ll_\eps \frac{1}{p^{\eps (M+1)/2}}\ll_M \frac{1}{p^{2}}.
$$

But for all  $\nub \in S^*(g)$ there exists 
$\betab \in {\cal B}(g)$ such that $\nub \geq \betab$. So the previous
relation implies that  the following bound is uniform  in $p$   and   $\s \in U_\eps $:
$$
\sum_{|\nub| \geq 1}\frac{f(p^{\nu_1},\dots,p^{\nu_n})}{p^{\la
    \sigmab ,\nub\ra}}
=\sum_{|\nub| \geq 1}\frac{g(\nub)}{p^{\la \sigmab ,\nub\ra}}
=\sum_{1\leq |\nub| \leq M}\frac{g(\nub)}{p^{\la \sigmab
    ,\nub\ra}} +O_\eps\left(\frac{1}{p^2}\right)
\ll_\eps \sum_{\nub \in {\cal B}(g) }
\frac{1}{p^{\la \sigmab ,\nub\ra}} 
+\frac{1}{p^2}.
$$
The multiplicativity of $f$ implies then that:
$$ 
{\cal M}(f;\s):=\sum_{\m \in \N^{n}}
\frac{f(m_1,\dots,m_n)}{m^{s_1}\dots m_n^{s_n}}=
\prod_p \left(\sum_{\nub \in \N_0^n} \frac{f(p^{\nu_1},\dots,p^{\nu_n})
    }{p^{\la \nub ,\s\ra}}\right)=
\prod_p \left(\sum_{\nub \in \N_0^n} \frac{g(\nub)}{p^{\la \nub
      ,\s\ra}}\right)
$$
converges absolutely in ${\cal D}_f \cap U_\eps$. \\
It follows that 
$s\mapsto {\cal M}(f;\s)$ converges absolutely in 
${\cal D}_f=\cup_{\eps >0} {\cal D}_f \cap U_\eps $.\par

$\bullet$ Let $\cb \in pol_0\left({\cal
    F}_0 (f)\right)$. The compactness of
the face ${\cal F}_0 (f)$ implies that $\cb \in \R_+^{*n}$. 
Moreover it follows from the definition of $pol_0\left({\cal
    F}_0 (f)\right)$ that 
$\forall \betab \in {\cal B}(g)$ $\la \betab , \cb\ra \geq 1$ with
equality if and only if $\betab \in {\cal F}_0 (f) \cap S^*(g)$.\\
Set $\delta_0 :=\frac{1}{2}\inf_{i=1,\dots,n} c_i >0$ and  
$\Omega_{\cb} :=\{\s \in \C^n \mid \forall i~ \Re(s_i) >c_i\}$.\\
Fix also $N:=\left[\frac{8}{\delta_0}+ \sup_{\x \in {\cal
      F}_0 (f)} |\x|\right]+1\in \N$. (Evidently, $N < \infty$ since ${\cal F}_0(f)$ is compact.) \\
It is easy to see that the following bound is uniform in   $p$   and  $\s \in
U_{-\delta_0} =\{\s \in \C^n \mid \forall i~ \Re(s_i) >-\delta_0\}$:\\
\begin{eqnarray*}
\sum_{|\nub| \geq N+1}\frac{g(\nub)}{p^{\la \cb+\sigmab ,\nub\ra}}
&\ll& \sum_{|\nub| \geq N+1}\frac{g(\nub)}{p^{\delta_0 |\nub|}}
\ll \frac{1}{p^{\delta_0 (N+1)/2}}
\sum_{|\nub| \geq N+1}\frac{g(\nub)}{p^{\delta_0 |\nub|/2}}\\
&\ll& \frac{1}{p^{\delta_0 (N+1)/2}}
\sum_{|\nub| \geq N+1}\frac{|\nub|^D}{2^{\delta_0 |\nub|/2}}
\ll \frac{1}{p^{\delta_0 (N+1)/2}}\ll \frac{1}{p^{2}}.
\end{eqnarray*}

Thus,  the following is  uniform  in $p$  and  $\s \in
U_{-\delta_0}$:
\begin{eqnarray}\label{remarkli}
\sum_{|\nub| \geq 1}\frac{f(p^{\nu_1},\dots,p^{\nu_n})}{p^{\la
    \cb+\s ,\nub\ra}}
&=&\sum_{|\nub| \geq 1}\frac{g(\nub)}{p^{\la \cb+\s ,\nub\ra}}
=\sum_{1\leq |\nub| \leq N}\frac{g(\nub)}{p^{\la \cb+\s
    ,\nub\ra}} +O\left(\frac{1}{p^2}\right) \nonumber \\
&=&\sum_{\nub \in {\cal F}_0(f)\cap S^*(g)}
\frac{g(\nub)}{p^{1+\la \s ,\nub\ra}} 
+\sum_{1\leq |\nub| \leq N \atop \nub \not 
\in {\cal F}_0(f)}
\frac{g(\nub)}{p^{\la \cb,\nub\ra +\la \s ,\nub\ra}} 
+O\left(\frac{1}{p^2}\right)
\end{eqnarray}
Since ${\cal F}_0(f)\cap S^*(g)$ is a finite set and  
$\la \cb, \nub \ra >1$ for all  
$\nub \in S^*(g) \setminus {\cal F}_0(f)$, 
it follows from (\ref{remarkli}) that there exists $\delta_1 \in ]0,\delta_0[$ and $\eps_1 >0$
such that the following is uniform in $p$  and   $\s \in
U_{-\delta_1}$:
\begin{equation}\label{eulerproddev}
\sum_{|\nub| \geq 1}\frac{f(p^{\nu_1},\dots,p^{\nu_n})}{p^{\la
    \cb+\s ,\nub\ra}}
=\sum_{|\nub| \geq 1}\frac{g(\nub)}{p^{\la \cb+\s ,\nub\ra}}
=\sum_{\nub \in {\cal F}_0(f) \cap S^*(g)}
\frac{g(\nub)}{p^{1+\la \s ,\nub\ra}} 
+O\left(\frac{1}{p^{1+\eps_1}}\right)
\end{equation}
The multiplicativity of $f$ now implies that ${\cal M}(f;\s)$
converges absolutely in\\
$\Omega_{\cb} :=\{\s \in \C^n \mid \forall i~ \Re(s_i) >c_i\}$ 
on which it can be written as follows: 
\begin{equation}\label{eulerprod=}
{\cal M}(f;\s):=\sum_{\m \in \N^{n}}
\frac{f(m_1,\dots,m_n)}{m^{s_1}\dots m_n^{s_n}}=
\prod_p \bigg(\sum_{\nub \in \N_0^n} \frac{f(p^{\nu_1},\dots,p^{\nu_n})
    }{p^{\la \nub ,\s\ra}}\bigg)=
\prod_p \bigg(\sum_{\nub \in \N_0^n} \frac{g(\nub)}{p^{\la \nub
      ,\s\ra}}\bigg)
\end{equation}
We conclude that for all $\s \in U_0 =\{\s \in \C^n \mid \forall i~
\Re(s_i) > 0\}$:
\begin{eqnarray}\label{ftfgfs}
G(f;\s)&:=&\bigg(\prod_{\nub \in {\cal F}_0(f) \cap S^*(g)} 
{\zeta (1+\la \nub , \s \ra )}^{-g(\nub)}\bigg) ~ {\cal
M}(f;\cb+\s)\nonumber \\
&=&\prod_p\bigg(\prod_{\nub \in {\cal F}_0(f) \cap S^*(g)} 
\left(1-\frac{1}{p^{1+\la \nub , \s \ra}}\right)^{g(\nub)}\bigg)~
\bigg(\sum_{\nub \in \N_0^n} \frac{g(\nub)}{p^{\la \nub ,\cb+\s\ra}}\bigg).
\end{eqnarray}
Combining (\ref{eulerproddev}) with (\ref{ftfgfs}) now implies that there
exist $\delta_2 \in ]0,\delta_1[$ and $\eps_2 >0$
such that :\\
$G(f;\s)=\prod_{p}\left(1+O\left(\frac{1}{p^{1+\eps_2}}\right)\right)$  
 uniformly in $\s \in U_{-\delta_2}$. 
It follows that the Euler product $\s \mapsto G(f;s)$ converges  
absolutely and defines a bounded holomorphic
function  on $U_{-\delta_2}$. \\
Moreover, since  
$\ds 
G(f;\zerob)=\prod_p
\bigg(1-\frac{1}{p}\bigg)^{\sum_{\nub \in {\cal F}_0 (f) \cap S^*(g) }g(\nub)} \cdot \bigg(\sum_{\nub \in \N_0^n} \frac{g(\nub)}{p^{\la \nub
      ,\cb\ra}}\bigg)$
is a convergent infinite product whose general term is $>0$, we conclude that $G(f;\zerob)>0$.\par
Moreover, for all $\s \in U_0 =\{\s \in \C^n \mid \forall i~
\Re(s_i) > 0\}$ we have:
\begin{eqnarray}\label{ftft0c}
H_{\cb}(f;\s)&:=& \left(\prod_{\nub \in {\cal F}_0 (f) \cap
    S^*(g)} {\la \nub , \s \ra }^{g(\nub)}\right) ~ {\cal M}(f;\cb+\s)
    \nonumber \\
&=& \left(\prod_{\nub \in  {\cal F}_0 (f) \cap
    S^*(g)}\left(\la \nub , \s \ra
    \zeta (1+\la \nub , \s \ra )\right)^{g(\nub)}\right)
    ~G(f;\s). 
\end{eqnarray}
Thus,  by using  the properties of the function $\s\mapsto G(f;\s)$
established above and standard   properties satisfied by  the Riemann zeta
function, it follows that there exists $\eps_0>0$ such that 
$\s \mapsto H_{\cb}(f;\s)$ has a holomorphic continuation with moderate
growth to $\{\s\in \C^n \mid \sigma_i >-\eps_0 ~\forall i\}$.\\
By  (\ref{ftft0c}) we conclude  that  $H_{\cb}(f;\zerob)= G(f;\zerob)>0$. 
This completes  the proof of   Proposition \ref{propofg}. \cqfd\par

{\bf Remark:} Let $\cb$ be a normalized polar vector of ${\cal F}_0(f)$. Assuming the compactness of  
${\cal F}_0(f)$   is equivalent to the assumption  $\cb
  \in (0,\infty)^n$. This condition 
is needed for our proof of Proposition \ref{propofg} above, and
  appears to be crucial for  the proof of Theorem A (see Appendix). 
Indeed, it does not yet seem possible to 
prove Theorem A without imposing this condition. \par

{\bf  Finishing the proof of Theorem \ref{ftfpourapplication}:}\\
Here we will use  the discussion given in 
the Appendix. Let $\cb \in pol_0\left({\cal F}_0(f)\right)$. As noted in the above  Remark, it follows that $\cb \in (0, \infty)^n.$  

Let ${\cal B}(g)$ be the set of $\betab \in S^*(g)$  which are in at
least one of the compact faces of the polyhedron ${\cal E}(f)$. Set  
$${\cal D}_f:= \{\s \in \C^n \mid \Re \la \s, \betab \ra >1 ~\forall \betab \in
{\cal B}(g)\}\cup \{ \s \in \C^n \mid \Re s_i > 0 ~\forall
i=1,\dots,n\} $$ 
and $$\Sigma_f:=\partial \left(\{\x \in \R^n \mid \la \betab, \Re \x \ra> 1 ~\forall
\betab \in {\cal B}(g)\} \right)\cap \R_+^{n}.$$ 

Proposition \ref{propofg} shows that ${\cal M}(f;s)$ converges in
${\cal D}_f$ and is, at worst, a  meromorphic function  at each point  
$\cb \in pol_0\left({\cal F}_0(f)\right)\cap \R_+^{*n}
\subset \Sigma_f \cap \R_+^{*n}$. It implies  also that the polar type of $f$ at $\cb$ is  
$\T_{\cb}=(I_{\cb}; \u ),$ where  
$I_{\cb}:= {\cal F}_0(f) \cap S^*(g)$ and   $u(\betab)= g(\betab)$ $\forall \betab
\in I_{\cb}$.\\ 
Let  
${\cal F}(\Sigma_f)(\unb)$ denote the face of $\Sigma_f$ with  polar vector
$\unb$. Lemma \ref{lienfaceetra} implies that $\iota (f) = \iota
({\cal E}(f)) =|\cb|$. We deduce that $|\cb|=\iota (\Sigma_f)$, and
consequently 
$\cb \in {\cal F}(\Sigma_f)(\unb) \cap \R_+^{*n}$. 

Moreover, by definition,   
${\cal F}_0(f)$ is the smallest
face of ${\cal E}(f)$ which meets the diagonal 
$\Delta =\R_+ \unb$. So, it is evident that 
$\unb  \in con^*(I_{\cb}).$  

In addition,  Proposition \ref{propofg} implies  that
$H_{\cb}(f;\zerob)\neq 0$.

Set now $r:=rank(I_{\cb})$ and fix in the sequel $\nub^1,\dots, \nub^r \in
I_{\cb}$ such that $rank\{\nub^1,\dots, \nub^r\}=r$.
Set also for any $\eps >0$ ${\cal R}(\eps) := \{\z =(z_1,\dots, z_r)\in \C^r
\mid |\Re (z_i) |< \eps ~\forall i\}$.\\
By using Theorem A (see Appendix),  the proof
of Theorem \ref{ftfpourapplication} will be complete once   we prove the
following:

{\it If $dim \,{\cal F}_0(f) = rank \left(S^*(g)\right)-1,$ 
    then there exist $\eps >0$ and a function $L$ such that:} 
\begin{equation}\label{Wholomorphe}
 L {\mbox {{\it is holomorphic in }}} {\cal R}(\eps) {\mbox {{\it and }}}   
H_{\cb}(f;\s)= L\left( \la \nub^1 , \s \ra,\dots, \la \nub^r, \s \ra
\right).
\end{equation}
Since we assume  that $dim \,{\cal F}_0 (f) = rank \left(S^*(g)\right)-1\,,$ 
it follows that $$rank(I_{\cb})= dim \,{\cal F}_0 (f)+1 =rank \left(S^*(g)\right).$$

By using relations  (\ref{ftft0c}) and (\ref{ftfgfs}), it is  clear
that it  suffices to prove that (\ref{Wholomorphe}) holds for the function 
$G(f;\s)$ (instead of $H_{\cb}(f;s)$), where 
$$
G(f;\s)=\prod_p\bigg(\prod_{\nub \in I_{\cb}} 
\left(1-\frac{1}{p^{1+\la \nub , \s \ra}}\right)^{g(\nub)}\bigg)~
\bigg(\sum_{\nub \in \N_0^n} \frac{g(\nub)}{p^{\la \nub ,\cb \ra+\la \nub ,\s\ra}}\bigg).
$$
First we recall from the proof of Proposition \ref{propofg} that there
exists $\delta_2>0$ such that the Euler product $G(f;\s)$ converges  
absolutely and defines a holomorphic
function  on $$U_{-\delta_2}:=\{\s \in \C^n \mid \Re(s_i) >-\delta_2 \ 
\forall i=1,\dots,n\}.$$  We fix this $\delta_2$ in the following discussion. 

It follows from the equality $rank\{\nub^1,\dots, \nub^r\}=r$ that the
  linear function $\varphi : \C^n \rightarrow \C^r,~ \s
\mapsto \varphi (\s)=(\la \nub^1, \s \ra,\dots,\la \nub^r, \s \ra )$
is onto. By a permutation of coordinates, if needed, we can then assume  that the function 
$\psi : \C^r \rightarrow \C^r,~ \s'=(s_1,\dots,s_r)
\mapsto \psi(\s'):=\varphi (\s,\zerob)=(\sum_{i=1}^r\nu_i^1
s_i,\dots,\sum_{i=1}^r \nu_i^r s_i)$ is an isomorphism.  
In particular there exist $\betab^1,\dots, \betab^r \in \Q^r$ so that 
$\forall \z \in \C^r,$\  $ \psi^{-1} (\z)=(\la \betab^1,\z\ra,\dots,
\la\betab^r,\z\ra)$, 
and 
\begin{equation}\label{psi-1}
\psi^{-1}\left({\cal R}(\eps_0)\right)\subset {\cal R}(\delta_2) {\mbox {
    for }} \eps_0 := \delta_2 \left(\max_{1\leq i\leq r}
    |\betab^i|\right)^{-1}>0.
\end{equation}
The fact that  $rank (S^*(g)) =rank (I_{\cb})= rank
\{\nub^1,\dots, \nub^r\}$ then  implies that for any $\nub \in S^*(g)$ there
exist $a_1(\nub),\dots, a_r(\nub) \in \Q$ such that 
$\nub = a_1(\nub) \nub^1+\dots +a_r(\nub)\nub^r$. It follows that for
all $\s \in U_{-\delta_2}$, $G(f;s)=W\left( \la \nub^1, \s \ra,\dots, \la \nub^r, \s \ra
\right)$ where  
\begin{equation}\label{G=W}
W(\z)=\prod_p\bigg(\prod_{\nub \in I_{\cb}} 
\left(1-\frac{1}{p^{1+\sum_{i=1}^r a_i(\nub)
      z_i}}\right)^{g(\nub)}\bigg)~ \bigg(1+\sum_{\nub \in S^*(g)} 
\frac{g(\nub)}{p^{\la \nub ,\cb \ra+\sum_{i=1}^r a_i(\nub)
    z_i}}\bigg).
\end{equation}
So to conclude, it suffices to prove that $W$ converges and defines a 
holomorphic function in ${\cal R}(\eps_0)$.  
But for all $\s' \in \psi^{-1}\left({\cal
    R}(\eps_0)\right)$, 
$${\tilde W} (\s')= W\circ \psi (\s')= W\left(\varphi
(\s',\zerob)\right)=G\left(f;(\s',\zerob)\right),$$ and by definition of
$\eps_0$ we have $\psi^{-1}\left({\cal R}(\eps_0)\right)\times
\{\zerob\}^{n-r} \subset {\cal R}(\delta_2) \times
\{\zerob\}^{n-r} \subset U_{-\delta_2}$. 
It follows that ${\tilde W}:= W\circ
\psi $ converges and defines a holomorphic function in 
$\psi^{-1}\left({\cal R}(\eps_0)\right)$. We deduce by composition that $W$ converges and defines a 
holomorphic function in ${\cal R}(\eps_0)$. This completes the proof of
Theorem \ref{ftfpourapplication}.\cqfd 
\subsection{Proofs of main results}
{\bf Proof of Theorems \ref{thMC} and \ref{thMCzeta}:}\\
By symmetry we have:
\begin{eqnarray*}
N_{H_P}(U({\bf A});t)&:=&\#\{M \in U({\bf A})\mid  {H_P} (M)\leq t\}
= c({\bf A}) \sum_{\{\m \in \N^{n+1}; P(\m)^{1/d} \leq t\}}
f(m_1,\dots, m_{n+1}).\\
Z_{H_P}(U({\bf A});s)&:=&\sum_{M \in U({\bf A})}  {H_P}^{-s}(M) 
= c({\bf A}) \sum_{\m \in \N^{n+1}} \frac{f(m_1,\dots, m_{n+1})}
{P(m_1,\dots, m_{n+1})^{s/d}}.
\end{eqnarray*}
where $c({\bf A})$ is the constant defined in (\ref{ca}),
and  $f$ is the function defined by:
\be
\item
$f(m_1,\dots,m_{n+1})=1$ if $m_1^{a_{i,1}}\dots m_{n+1}^{a_{i,n+1}}=1$ 
$\forall i=1,\dots , l$ and gcd$(m_1,\dots,m_{n+1})=1$  
\item 
$f(m_1,\dots,m_{n+1})=0$ otherwise.
\ee
It is easy to see that $f$ is a multiplicative function and that for
any prime number $p$ and any $\nub \in \N_0^{n+1}$: 
$f(p^{\nu_1},\dots ,p^{\nu_{n+1}})= g(\nub)$ where $g$ is the
characteristic function of the set $T({\bf A})$ defined in \S 3. 
So, it is obvious that $f$ is a  uniform multiplicative function.\\ 
It now suffices to verify that the assumptions of Theorem
\ref{ftfpourapplication} and Corollary \ref{applicationmulti} from \S 4.1 are satisfied.\\ 
By using the notations of \S 4.1, it is easy to check that  
$${\cal E}(f)= {\cal E}\left(T^*(\bold A)\right) {\mbox { and }} {\cal F}_0(\bold
A)={\cal F}_0 (f).$$ 
Therefore, the ellipticity of $P$ and the compactness 
of the face ${\cal F}_0(\bold A)$  imply that the first part of
Theorem \ref{thMCzeta} (resp. Theorem \ref{thMC}) follows  from
Theorem \ref{ftfpourapplication} (resp. Corollary \ref{applicationmulti}).\par
Let us now suppose that $dim\left({\cal F}_0(\bold A)\right)=dim X(\bold
A)=n-l$. Since  
$$dim\left({\cal F}_0(\bold A)\right)=rank\left({\cal
    F}_0(\bold A)\cap T^*(\bold A)\right) -1 \leq 
rank\left(T^*(\bold A)\right)-1\leq dim \R^{n+1}-rank(\bold A)-1=n-l,$$ it follows that  
 $dim\left({\cal F}_0(\bold A)\right)=rank\left(T^*(\bold A)\right)-1
=rank\left(S^*(g)\right)-1$. Consequently, the second part of Theorem
\ref{thMCzeta} (resp. Theorem \ref{thMC}) also 
follows from  Theorem \ref{ftfpourapplication} (resp. Corollary \ref{applicationmulti}).\cqfd \par
{\bf Proof of Theorem \ref{CMa}:}\\
Let $A_n(\a)$ be the $1\times (n+1)$ matrix $A_n(\a):=(a_1,\dots , a_n,-q)$.
It is then clear that   $X_n(\a)=X\left(A_n(\a)\right).$ Thus,   
Theorem \ref{CMa} will follow directly   from   Corollary
\ref{applicationmulti} once we show that all the hypotheses of the corollary are satisfied.\\
Set $c(\a):= \frac{1}{2}\#\left\{(\eps_1,\dots ,\eps_{n+1})\in \{-1,+1\}^{n+1}
\mid \eps_1^{a_1}\dots \eps_n^{a_n}=\eps_{n+1}^q \right\}$.\\ 
  As above,   by symmetry we have for all $t>1$:
$$N_{H_{P}}(U_n(\a);t)=\#\{M \in U_n(\a) \mid  H_{P}(M)\leq t\}
= c(\a) \sum_{\{\m \in \N^{n}; {\tilde P}(\m)^{1/d} \leq t\}}
f(m_1,\dots, m_{n})
$$
where ${\tilde P}$ is the generalized polynomial defined by 
${\tilde P}(X_1,\dots,X_n):=P\left(X_1,\dots,X_n, \prod_j X_j^{a_j/q} \right)$ and $f$ is the function defined by:
\be
\item
$f(m_1,\dots,m_{n})=1$ if $m_1^{a_1}\dots m_n^{a_n} $ is the $q^{th}$
power of an integer and gcd$(m_1,\dots,m_{n})=1$;  
\item 
$f(m_1,\dots,m_{n})=0$ otherwise.
\ee
It is easy to see that $f$ is a multiplicative function and that for
any prime number  $p$ and any $\nub \in \N_0^{n}$: 
$f(p^{\nu_1},\dots ,p^{\nu_{n}})= g(\nub)\,,$ where $g$ is the
characteristic function of the set $L_n(\a) \cup \{\zerob\}$ and     
$$L_n(\a):= \left\{\r \in 
\N_0^n \setminus \{\zerob\};~  q | \la \a, \nub \ra {\mbox { and }} \nu_1\dots \nu_n =0 \right\}.$$ 
So it is clear that   $f$ is also uniform.\\
By definition, we have that 
${\cal E}(f) = {\cal E} (\a) = {\cal E}\left(L_n(\a)\right).$ 
A simple check verifies that 
$rank\left(S^*(g)\right)= rank\left(L_n(\a)\right)=n.$ Moreover,  
the face ${\cal F}_0(\a) = {\cal F}_0(f)$ is compact.  
Consequently the assumption $dim {\cal
  F}_0 (f) = rank\left(S^*(g)\right)-1$ 
is equivalent to the assumption that ${\cal F}_0 (\a)$  is a facet of 
${\cal E}(\a)$. As a result, Theorem \ref{CMa} follows from Corollary
\ref{applicationmulti}, once one has also noted  that   Lemma
\ref{lienfaceetra} implies   
$\iota (f)=\iota \left({\cal E}(f)\right)=|\bold c|$ for any normalized polar vector $\bold c$ of ${\cal F}_0 (\a).$ \cqfd \par

{\bf Remark 4:} The two examples of a uniform multiplicative function used above are, evidently, quite special. The reader may wonder if other examples exist that are not merely characteristic functions of subsets of lattice points. One such example is found in the article by Swinnerton-Dyer  \cite{swd}, in which he was interested in finding a ``canonical" height function for the surface $X_1 X_2 X_3 = X_4^3.$ 
His idea  was to define a height function at each finite prime
  of the form 
$$f (p^{\nu_1},p^{\nu_2}, p^{\nu_3}, p^{\nu_4}) = g(\nu_1, \nu_2, \nu_3, \nu_4) \cdot {\mbox {characteristic function of }} \{\nu_1 + \nu_2 + \nu_3 = 3 \nu_4\},$$
where $g$ was a quadratic polynomial in $(\nu_1 - \nu_4, \nu_3 - \nu_4, \nu_3 - \nu_4).$ It follows that this $f$ is a uniform multiplicative function. 

The height   zeta function for which   Swinnerton-Dyer was able to
prove the asymptotic density was determined  by an adelic metric $v = (v_p)$ for $\P^3 (\Q)$ with $v_\infty = 0.$ 
It would be quite interesting to introduce the {\it nonzero} metric
for $v_\infty$ that he  also proposed, and derive the asymptotic for
the ``canonical" height density. Swinnerton-Dyer did not know how to
do this. Perhaps the methods of this paper, suitably extended,  can be used to do this.  

\section*{Appendix on   mixed zeta functions}
In order to make this paper self contained, we state, for the
convenience of the reader, the main result   of  \cite{essouamixte},
from which Theorem \ref{ftfpourapplication} follows. We do so in a somewhat simplified  
context that suffices to prove Theorem \ref{ftfpourapplication}.\par
Let $f:\N^{n}\rightarrow \C$ be any function. Suppose that there exists a
finite subset $I$ of $\R_+^n\setminus
  \{\zerob\}$  such that ${\cal M}(f;\s):=\sum_{m_1,\dots,m_n \geq 1}
\frac{f(m_1,\dots,m_n)}{m_1^{s_1}\dots m_n^{s_n}}$ converges absolutely in  
$${\cal D}_f:=\{\s \in \C^n \mid \la \betab, \Re \s \ra> 1 ~\forall
\betab \in I\}\cap \{\s \in \C^n \mid \Re(s_i)>0 ~\forall
i=1,\dots,n\}.$$
Set $\Sigma_f:=\partial \left(\{\x \in \R^n \mid \la \betab, \Re \x \ra> 1 ~\forall
\betab \in I\}\right)\cap \R_+^{n}\subset \partial {\cal D}_f$ and   
${\cal F}(\Sigma_f)(\unb)$ the face with polar vector $\unb$ 
(see \S 2.2).\par

{\bf Definition A.}
{\it The function $f:\N^{n}\rightarrow \C$ is said to be of
  finite type if there exists a point 
$\cb \in {\cal F}(\Sigma_f)(\unb)$  such that $\ s\mapsto {\cal M}(f;\s)$ 
has a meromorphic continuation with moderate growth to a neighbourhood
of $\cb$ in the following sense:\\
Let  $I_{\cb}:=\{\betab  \in I \mid \la \betab ,\cb \ra =1\}.$  There
exists a family $u = \big(u(\betab)\big)_{\betab \in
  I_{\cb}}$ of positive integers and $\eps_0 >0$ such that:
$\ds \s\mapsto H_{\cb}(f; \s):=\big(\prod_{\betab \in I_{\cb}} {\left(\la
    \s,\betab\ra \right)}^{u(\betab)}\big)~{\cal M}(f;\cb+\s)$ 
has a {\it holomorphic} continuation with moderate growth in $\Im \s$ to the set  
$\{ \s \in \C^n \mid \sigma_i >-\eps_0 ~\forall i\}$, and does not
  vanish identically along any hyperplane $\la \s,\betab\ra = 0$ $(\betab \in
I_{\cb})$.\par
We call the pair $\T_{\cb} :=(I_{\cb},\u ),$ the polar type of $f$
at $\cb$.}\par 

{\bf Main result on the mixed zeta functions (simplified version):}\\
Let $P=P(X_1,\dots, X_n)$ be a (generalized)  polynomial 
with positive coefficients that is elliptic and  homogeneous of degree
$d>0$. Let $f$ be a function of finite type. 
Set \\  
$\ds Z(f;P;s):=\sum_{\m \in \N^{n}}
\frac{f(m_1,\dots,m_n)}{P(m_1,\dots,m_n)^{s/d}}$.\\ 
 Theorems 1 and 2, proved in \cite{essouamixte}, imply the following: \par 
{\bf Theorem A.}
{\it Let $\cb \in {\cal F}(\Sigma_f)(\unb)$ be a point at which the property in Definition A is satisfied by ${\cal M}(f;\s)$ in some neighborhood of $\cb.$ \   
Then, $s \mapsto Z(f;P;s)$ is a holomorphic function in the half-plane 
$\{s: \sigma > |\cb|\},$ and there exists $ \eta >0$ such that   
 $s \mapsto Z(f;P;s)$ has a meromorphic continuation with moderate
 growth to the half-plane $\{\sigma > |\cb| - \eta\}$ with at most
one pole at $s=|\cb|$ of order at most equal to $\rho_0:= \sum_{\betab \in I_{\cb}} u(\betab)
-rank(I_{\cb}) +1$.\par
Suppose in addition that:
\be
\item $\unb \in con^*(I_{\cb})$; 
\item there exists a function $L$ holomorphic in  neighborhood of 
  $\zerob$ such that:\\ 
$\ds {\cal M}(f;\s)= L\left((\la \betab,\s\ra)_{\betab \in
 I_{\cb}}\right);$ 
\item $H_{\cb}(f;\bold 0) \neq 0.$
\ee
Then  
$s=|\cb|$ is indeed a pole of order $\rho_0$ and   
$$Z(f;P;s)\sim_{s\rightarrow |\cb|}
\frac{C(f;P)}{\left(s-|\cb|\right)^{\rho_0 }}\,, \qquad {\mbox { where
  }}\quad 
C(f;P):= H_{\cb}(f;\zerob) d^{\rho_0 } A_0(\T_{\cb}, P)\neq 0.$$ 
($A_0(\T_{\cb}, P)$ is the mixed volume constant (see \S 2.3.3)
associated to $\T_{\cb}$ and $P$).
}\par 
{\bf Acknowledgments:} The author wishes to express his thanks to Ben
Lichtin for his many helpful suggestions and his careful reading of this paper.


\end{document}